\documentclass[12pt]{article}

\usepackage{graphicx}
\usepackage{amsmath}
\usepackage{amsfonts}
\usepackage{amssymb}
\usepackage{ifthen}

\newtheorem{thm}{Theorem}[section]
\newtheorem{lemma}[thm]{Lemma}
\newtheorem{defn}[thm]{Definition}
\newtheorem{defns}[thm]{Definitions}

\newtheorem{example}[thm]{Example}

\newtheorem{propn}[thm]{Proposition}
\newtheorem{cor}[thm]{Corollary}

\newcommand{\os}{\le_O}
\newcommand{\ons}{\triangleleft_O}
\newcommand{\cs}{\le_C}
\newcommand{\cns}{\triangleleft_C}

\newcommand{\nons}{{\not \triangleleft}_O}

\newcommand{\Int}{\mathop{Int}}

\newcommand{\T}{\mathop{Tor}}
\newcommand{\Hom}{\mathop{Hom}}

\newcommand{\Cr}{\prod}
\newcommand{\Dr}{\bigoplus}

\newcommand{\cen}{\mathop{C}}
\newcommand{\zg}{\mathop{Z}}

\newcommand{\norm}{\mathop{N}}
\newcommand{\Aut}{\mathop{Aut}}

\newcommand{\sub}{\mathcal{S}}
\newcommand{\nsub}{\mathcal{N}}

\newcommand{\h}{\mathop{h}}

\newcommand{\df}[1]{\textbf{#1}}

\newcommand{\hht}{\mathop{ht}}

\newcommand{\Der}{\mathop{Der}}

\newenvironment{prf}[1][Proof]{\textbf{#1.}
}{\ \rule{0.5em}{0.5em} \vskip .5em}
\begin{document}

\title{Counting the Closed  Subgroups \\
 \hfill of Profinite Groups}



\date{November 2007}
\author{Paul Gartside and Michael Smith}

\maketitle

\begin{abstract}
The  sets of closed and closed-normal subgroups of a
profinite group carry a natural profinite topology. Through a
combination of algebraic and topological methods the size of these
subgroup spaces is calculated, and the spaces partially classified up
to homeomorphism.  
\end{abstract}

\section{Introduction}
In this paper we calculate the possible cardinalities of $\sub (G)$,
the set of all closed subgroups of a profinite group $G$, and
find conditions on $G$ determining the cardinality of $\sub (G)$.
In summary: $\sub (G)$ is finite if and only if $G$ is finite; $\sub
(G)$ is countably infinite if and only if $G$ is a finite central
extension of $\Dr_{i=1}^n\mathbb Z_{p_i}$ where the  $p_i$ are
distinct primes and $\mathbb{Z}_{p_i}$ is a copy of the $p_i$-adic integers; and
otherwise $|\sub (G)| = 2^{w(G)}$ where $w(G)$ is the weight of $G$
(the cardinality of a minimal sized base for the topology of
$G$). 

We further show that for a profinite group $G$, the set $\nsub (G)$ of
closed normal subgroups of $G$ is either countable or of size
$2^{w(G)}$.

These results follow from a mixture of algebraic and topological
considerations. In particular we use a natural topology defined on $\sub (G)$,
and inherited by $\nsub (G)$, introduced in \cite{FG1} and explored
further in \cite{FG2}. Extending our results on the number of
closed (normal) subgroups we also consider the problem of classifying
the space of closed subgroups $\sub (G)$ and closed normal subgroups
$\nsub (G)$ up to homeomorphism. Here we give a complete solution in
the case when $\sub (G)$ is countable and when $w(G) = \aleph_1$ (the
first uncountable cardinal -- in this case $\sub (G)$ is homeomorphic
to $\{0,1\}^{\aleph_1}$). The situation when $w(G) > \aleph_1$ is
rather mysterious: in particular $\sub (G)$ need not be homeomorphic to
$\{0,1\}^{w(G)}$. We also determine when $\sub (G)$ is homeomorphic to the
Cantor set (equivalently, homeomorphic to $\{0,1\}^{\aleph_0}$). These
results all depend on determining when  a subgroup 
space has an isolated point.

The remaining case in the topological classification of the space of
closed subgroups of a profinite group is the `mixed' case: the
subgroup space is uncountable, has a countable base, but is not
homeomorphic to the Cantor set. This case is investigated further in
\cite{GS2}.

\section{Background Material}

\paragraph{Subgroup Spaces}
 By
definition, a {\sl profinite} group $G$ is one that can be represented as
a {\sl projective} limit of {\sl finite} groups, $\lim_{\leftarrow}
G_{\lambda}$. Writing $\sub (G)$ 
for the set of closed subgroups of a profinite group,
$G=\lim_{\leftarrow} G_{\lambda}$, one sees that $\sub (G) =
\lim_{\leftarrow} \sub (G_\lambda)$. Giving the finite set
$\sub(G_\lambda)$ the discrete topology for each $\lambda$, we see
that the projective limit $\sub(G)$ picks up a natural topology. This
topology is profinite (compact, Hausdorff and zero-dimensional).

An alternative description of the topology on 
$\sub (G)$ is that it is the subspace topology inherited by $\sub(G)$ from the
space of all compact subsets of $G$ with the Vietoris topology, and so
the topology is independent of the particular projective
representation of $G$.

We can concretely describe canonical basic open neighbourhoods of a subgroup
in $\sub (G)$ or $\nsub (G)$ for a profinite group $G$ as follows.
\begin{defn}
Let $G$ be a profinite group. For $H\cs G$ and $N\ons G$, we define
$B(H,N)=\{K\cs G\mid KN=HN\}$.
\end{defn}
\begin{lemma}\label{lclbase}
Let $G$ be a profinite group, and $H\cs G$. Suppose that
$(N_{\lambda})_{\lambda\in\Lambda}$ is a family of open normal subgroups of
$G$, forming a base for the open neighbourhoods of $1$ in $G$. Then
$(B(H,N_{\lambda}))_{\lambda\in\Lambda}$ forms a base for the open
neighbourhoods of $H$ in $\sub (G)$. Also, if $H\cns G$ then
$(B(H,N_{\lambda})\cap\nsub (G))_{\lambda\in\Lambda}$ forms a base for the
open neighbourhoods of $H$ in $\nsub (G)$.
\end{lemma}
\begin{lemma}\label{lsgchwt}
Let $G$ be an infinite profinite group. Then
$\chi (1,\sub (G))=\chi (1,\nsub (G))=w(G)$.
\end{lemma}
\begin{cor}\label{csubwt}
Let $G$ be an infinite profinite group. Then $w(\sub (G))=w(\nsub (G))=w(G)$.
\end{cor}

\begin{propn}\label{pctprfct}
If a profinite group $G$ is topologically isomorphic to $\Cr_{i=1}^n
G_{p_i}$ where $p_1, \ldots , p_n$ are distinct primes, and $G_{p_i}$
is a pro-$p_i$ group, then $\sub (G) \cong \prod_{i=1}^n \sub(G_{p_i})$.
\end{propn}

\paragraph{Classifying Profinite Spaces}

\begin{defn}
Let $X$ be a topological space. For $Y\subseteq X$, let $Y^{\prime}$ denote
the set of all limit points of $Y$, that is $Y\setminus Y^{\prime}$ is the
set of points of $Y$ which are isolated in $Y$. We define the following
transfinite sequence.

$X^{(0)}=X$,

$X^{(\alpha +1)}=(X^{(\alpha)})^{\prime}$, for $\alpha$ an ordinal,

$X^{(\lambda)}=\bigcap_{\mu <\lambda}X^{(\mu)}$, for $\lambda$ a limit ordinal.
\end{defn}
We now list some facts about this \df{Cantor--Bendixson} process.
\begin{lemma}\label{lctbdx}
Let $X$ be a Hausdorff space.
\begin{enumerate}
\item[(i)] $X^{(\alpha)}$ is closed in $X$ for every ordinal $\alpha$.
\item[(ii)] If $\alpha\leqslant\beta$ then $X^{(\alpha)}\supseteq X^{(\beta)}$.
\item[(iii)] $X^{(\alpha)}\setminus X^{(\alpha +1)}$ is the set of isolated
points of $X^{(\alpha)}$ for every ordinal $\alpha$, and is countable if
$X$ is countably based.
\item[(iv)] $(X^{(\alpha)})^{(\beta)}=X^{(\alpha +\beta)}$ for all ordinals
$\alpha$ and $\beta$.
\item[(v)] If $Y\subseteq Z\subseteq X$ then
$Y^{(\alpha)}\subseteq Z^{(\alpha)}$ for every ordinal $\alpha$.
\item[(vi)] There is a least ordinal $\lambda$ such that
$X^{(\lambda)}=X^{(\lambda +1)}$ and if $\alpha\geqslant\lambda$ then
$X^{(\alpha)}=X^{(\lambda)}$. $X^{(\lambda)}$ is perfect in itself. If $X$ is
countably based then $\lambda$ is countable.
\item[(vii)] If $X$ is compact, and $X^{(\alpha)}=\emptyset$ for some ordinal
$\alpha$, then there is a successor ordinal $\lambda$ such that
$X^{(\lambda)}=\emptyset$ and $X^{(\lambda -1)}$ is a finite non-empty
discrete space.
\item[(ix)] If $Y$ is open in $X$ then $Y^{(\alpha)}=Y\cap X^{(\alpha)}$ for
every ordinal $\alpha$.
\end{enumerate}
\end{lemma}
It follows from Lemma \ref{lctbdx} that every countably based Hausdorff
space can be written as the disjoint union of a countable set (the
\df{scattered part} of $X$) and a set which is perfect in itself (the
\df{perfect hull} of $X$); this is the Cantor--Bendixon theorem.
\begin{propn}\label{pchcant}
Let $X$ be a non-empty profinite space. Then $X$ is homeomorphic to the
Cantor set if and only if $X$ is perfect and countably based.
\end{propn}

\begin{defns}
Let $X$ be a Hausdorff space. The \df{scattered height} of $X$, $\hht (X)$ is
the least ordinal $\lambda$ such that $X^{(\lambda)}=X^{(\lambda +1)}$.

Let $x\in X\setminus X^{(\hht (X))}$. $\hht (x,X)$ is defined to be the least
ordinal $\alpha$ such that $x\not\in X^{(\alpha +1)}$.
\end{defns}
Of course to say that $x\in X\setminus X^{(\hht (X))}$ has $\hht (X)=\alpha$
is precisely the same as saying that $x\in X^{(\alpha)}$ and that $x$ is
isolated in $X^{(\alpha)}$. The next lemma follows immediately from the
definitions.

\begin{lemma}\label{lhtopen}
Let $X$ be a Hausdorff space, $Y$ an open set in $X$ and $x\in Y$. Then
$\hht (Y)\leqslant \hht (X)$ and $\hht (x,X)=\hht (x,Y)$.
\end{lemma}

In relation  to Lemma \ref{lctbdx}(vii), a compact Hausdorff space $X$ for
which $X^{(\hht (X))}=\emptyset$ is sometimes known as a \df{scattered} space.
It is clear that a countably based profinite space is scattered if and only
if it is countable.

\begin{lemma}
Let $\alpha$ be a non-zero countable ordinal, and $n$ be a positive integer.
Then $\omega^{\alpha}n+1$, with the order topology,  is a countably
infinite profinite space. Moreover 
$(\omega^{\alpha}n+1)^{(\alpha)}=\{\omega^{\alpha},\omega^{\alpha}2,\ldots
,\omega^{\alpha}n\}$. In particular $\hht (\omega^{\alpha}n+1)=\alpha +1$
and $|(\omega^{\alpha}n+1)^{(\alpha)}|=n$.

Note that $\omega+1$ is homeomorphic to a convergent sequence.
\end{lemma}

\begin{propn}\label{pctprs}
Let $X$ be a countable profinite space. Then $X$ is homeomorphic to
$\omega^{\hht (X)-1}|X^{(\hht (X)-1)}|+1$.
\end{propn}

\section{When $\sub (G)$ is Countable}

In this section we characterise the profinite groups which have precisely
countably infinitely many closed subgroups. {\sl A fortiori} such a
group has countably many open subgroups, so  must be countably based.
The following theorem (originally established by Dikranjan \cite{dik}
and re--discovered by Morris, Oates-Williams and Thompson
\cite{MOWT}) and example suggests
that these groups might be built up in some simple way from the
$p$-adic integers $\mathbb Z_p$. Our
main result, Theorem \ref{tcount}, confirms this.
\begin{thm}\label{tmowt}
Let $G$ be an infinite compact Hausdorff topological group. Then the following
are equivalent.
\begin{enumerate}
\item[(i)] Every non-trivial closed subgroup of $G$ is open.
\item[(ii)] $G$ is topologically isomorphic to $\mathbb Z_p$ for some prime
$p$.
\end{enumerate}
\end{thm}
\begin{example}[See \cite{FG1}]\label{esubzp}
$\sub (\mathbb{Z}_p) \cong (\omega +1)$.

More precisely, the open subgroups $p^n \mathbb{Z}_p$ form a sequence
converging to the sole non-open closed subgroup, the trivial subgroup.
\end{example}

For infinite profinite groups, having countably infinitely many closed
subgroups turns out to be equivalent to having less than $2^{\aleph_0}$
closed subgroups. To show this we do not need to assume the Continuum
Hypothesis.  

This fact will
follow immediately from the topological considerations of 
Section~\ref{sspcancub} (see Corollary~\ref{ccount}).
Here we give a direct group-theoretic argument. 

\paragraph{Continuum Many Closed Subgroups} 
First we give some conditions which ensure that a profinite group has
at least $2^{\aleph_0}$ closed subgroups.
\begin{lemma}\label{lexsm}
\mbox{\ }
\begin{enumerate}
\item[(i)] Let $(G_i)_{i\in I}$ be an infinite family of profinite groups where
infinitely many of the $G_i$ are non-trivial. Then $G=\Cr_{i\in I}G_i$ has at
least $2^{\aleph_0}$ closed subgroups.
\item[(ii)] Let $G$ be an infinite abelian torsion profinite group. Then $G$
has at least $2^{\aleph_0}$ closed subgroups.
\item[(iii)] Let $G$ be a profinite group, $N\cns G$, and $K\cns G$ with
$K\leqslant N$. Suppose for some prime $p$, $K\cong\mathbb Z_p$ and
$G/N\cong\mathbb Z_p$. Then $G$ has at least $2^{\aleph_0}$ closed subgroups.
\end{enumerate}
\end{lemma}

\begin{prf}
For (i), let $J$ be a countably infinite subset of $I$ such that for every
$i\in J$, $G_i$ is non-trivial. Then for each of the
$2^{\aleph_0}$-many distinct subsets $J'$ of $J$, we have distinct
closed subgroups $G_{J'} = \Cr_{i \in I} G_i(J')$ where $G_i(J')=G_i$
if $i \in J'$ and $G_i(J')=1$ otherwise.

(ii) follows immediately from (i) and the fact that compact abelian
torsion groups are topologically isomorphic to a Cartesian sum of
finite cyclic groups of bounded order.

Now, for (iii), let $(G/K)/(N/K)=\overline{\langle (N/K)Kx\rangle}$,
$L/K=\overline{\langle Kx\rangle}$, and let $P/K$ be a Sylow $p$-subgroup
of $L/K$. Since
$(L/K)/(N/K\cap L/K)\cong (N/K)(L/K)/(N/K)=(G/K)/(N/K)\cong G/N\cong\mathbb Z_p$
, $p^{\infty}\mid o(L/K)$. Since $L/K$ is procyclic and $P/K$ is not
finite, we see that $P/K\cong\mathbb Z_p$. 
Hence without loss of generality we may assume that $K=N$ (and
thus $P=G$). We now show that $G$ splits over $N$. Let
$G/N=\overline{\langle Na\rangle}$. Then $G=N\overline{\langle a\rangle}$ and
so $\overline{\langle a\rangle}/\overline{\langle a\rangle}\cap N\cong N
\overline{\langle a\rangle}/N=G/N\cong\mathbb Z_p$. So
$\overline{\langle a\rangle}\cap N=1$ and $G=N\rtimes\overline{\langle a\rangle}
$.

Now define $\phi\colon N\to\{H\mid H\cs G\}$, by
$\phi (n)=\overline{\langle na\rangle}$. Suppose $n_1\neq n_2$ but
$\phi (n_1)=\phi (n_2)$. Then as
$(n_1 a)(n_2 a)^{-1}\in\overline{\langle n_1 a\rangle}$ and
$(n_1 a)(n_2 a)^{-1}=n_1n_2^{-1}\in N\setminus\{1\}$,
$\overline{\langle n_1 a\rangle}\cap N\neq 1$. As
$\overline{\langle a\rangle}\cap N=1$, $n_1 a\not\in N$, and so
$\overline{\langle n_1 a\rangle}\nleqslant N$. Thus
$\overline{\langle n_1 a\rangle}N/N$ is a non-trivial subgroup of $G/N$, and
as $\mathbb Z_p\cong\overline{\langle n_1 a\rangle}N/N\cong\overline{\langle
n_1 a\rangle}/(\overline{\langle n_1 a\rangle}\cap N)$,
$\overline{\langle n_1 a\rangle}\cap N=1$, a contradiction. So $\phi$ is
injective, and since $|N|=2^{\aleph_0}$, $G$ has at least $2^{\aleph_0}$ closed
subgroups.
\end{prf}

So if a profinite group has less than $2^{\aleph_0}$ closed subgroups it
cannot have a closed section of any of the above three forms.

\paragraph{The Abelian Case and the Center}
Next we shall prove a special case of part of the main result; namely, a
necessary condition for an infinite abelian profinite group to have less than
$2^{\aleph_0}$ closed subgroups. This result can be obtained using
Pontryagin duality from a result of Boyer (Lemma 3 of \cite{Boyer}) or
from a result of Berhanu, Comfort and Reid (Proposition 2.3 of \cite{BCR}).
Here we give a direct argument.
\begin{lemma}\label{lcab}
Let $G$ be an infinite abelian profinite group. Suppose $G$ has less than
$2^{\aleph_0}$ closed subgroups. Then there is a finite non-empty set of
primes $\{p_1,\ldots ,p_n\}$ (with $p_i\neq p_j$ for $i\neq j$) and a finite
abelian group $T$ such that $G$ is topologically isomorphic to
$T\times \Dr_{i=1}^n \mathbb Z_{p_i}$.
\end{lemma}
\begin{prf}
Since $G$ is abelian it is topologically isomorphic to the
Cartesian product of its Sylow subgroups. 
So by Lemma \ref{lexsm}(i),
there are only finitely many primes, $p_1, \ldots, p_n$, dividing the
order of $G$, and  
$G = \Cr_{i=1}^n G_{p_i}$, where $G_{p_i}$ is pro-$p_i$. Hence
$\sub(G) = \prod_{i=1}^n \sub (G_{p_i})$, so it now suffices to prove
  the result in the case 
where $G$ is a pro-$p$ group for some prime $p$.

By Lemma \ref{lexsm}(ii) $G$ is not a torsion group. Let $x$ be an element
of $G$ of infinite order, and let $H=\overline{\langle
  x\rangle}$. Then
 $H\cong\mathbb Z_p$. Let $y\in G$, and let
$K=\overline{\langle x,y\rangle}$. Then $K$ is a finitely generated abelian
pro-$p$ group, so
$K=A\times B$ where $A$ is a finitely generated torsion-free abelian pro-$p$
group and $B$ is a finite abelian $p$-group. $A\cong\mathbb Z_p^m$ for some
positive integer $m$. But by Lemma \ref{lexsm}(iii), $m=1$, and so $A=H$.
Let $T=\T (G)$. Then $B\leqslant T$, and so $K=A\times B\leqslant H\times T$.
Thus for any $y\in G$, $y\in H\times T$. Hence $G=H\times T$.

We now show that $T$ is finite. For each non-negative integer $i$, let
$T_i=\{g\in G\mid g^{p^i}=1\}$. Then $T_i\cs G$ as $T_i$ is the kernel of the
continuous endomorphism $g\mapsto g^{p^i}$. As
$T=\bigcup_{i\geqslant 0}T_i$, $G=\bigcup_{i\geqslant 0}(H\times T_i)$, and
so by Baire's Category theorem there is an $i$ such that $H\times T_i\os G$.
Thus there is a non-negative integer $k$ such that $|G:H\times T_i|=p^k$.
Now for all integers $j\geqslant i$,
$p^k\geqslant |HT_j:HT_i|=|T_j:(HT_i)\cap T_j|=|T_j:HT_i|$. It follows
that $T=T_{i+k}$. Thus $T\cs G$. By Lemma \ref{lexsm}(ii), $T$ is finite.
\end{prf}

\begin{lemma}\label{labvp}
Let $G$ be an abelian profinite group and $N\os G$. Suppose for some finite
non-empty set of primes $\{p_1,\ldots ,p_n\}$ (with
$p_i\neq p_j$ for $i\neq j$) that $N\cong\Dr_{i=1}^n \mathbb Z_{p_i}$. Then
there is a finite abelian group $T$ such that $G$ is topologically isomorphic
to $T\times\Dr_{i=1}^n \mathbb Z_{p_i}$.
\end{lemma}

\begin{prf}
Again as $G$ is topologically isomorphic to the Cartesian product of its Sylow
subgroups, it suffices to prove the result in the case where $G$ is a
pro-$p$ group for some prime $p$. Let $N=\overline{\langle x\rangle}$. As in
the proof of Lemma \ref{lcab} above, let $y\in G$, and let
$K=\overline{\langle x,y\rangle}$. Again $K=A\times B$ where $A$ is a finitely
generated torsion-free abelian pro-$p$ group and $B$ is a finite abelian
$p$-group. Since $\h (G)=1$, $A=N$ and as above $K\leqslant N\times\T (G)$. So
$G=N\times\T (G)$.
\end{prf}

The next proposition is well known.

\begin{propn}\label{pzgocomf}
\mbox{\ }
\begin{enumerate}
\item[(i)] If $G$ is a group with $|G:\zg (G)|<\infty $, then $G^{\prime}$
is finite.
\item[(ii)] If $G$ is a finitely generated profinite group with $G^{\prime}$
finite, then $\zg (G)\ons G$.
\end{enumerate}
\end{propn}

\begin{prf}
(i) is a well-known result of Schur; see, for example 10.1.4 of \cite{Rob}.

For (ii), suppose that $G=\overline{\langle g_1,\ldots ,g_n\rangle}$. Now for
each $i$, with $1\leqslant i\leqslant n$, define
$\phi_i\colon \cen_G(G^{\prime})\to G$ by $\phi_i(x)=[x,g_i]$ for
$x\in\cen_G(G^{\prime})$. Now each $\phi_i$ is a continuous homomorphism since
$\phi_i(x_1x_2)=[x_1,g]^{x_2}[x_2,g]=[x_1,g][x_2,g]$, for
$x_1,x_2\in\cen_G(G^{\prime})$. Since $G^{\prime}$ is a finite normal subgroup
of $G$, $\cen_G(G^{\prime})\ons G$. Hence, for each $i$, $\ker (\phi_i)\os G$.
Let $K=\bigcap_{i=1}^n\ker (\phi_i)$. Then $K\os G$. For each $i$,
$g_i\in\cen_G(K)$. Thus $K$ is central in $G$, and so $\zg (G)\ons G$ as
required.
\end{prf}

\paragraph{General Case}
We are now ready to characterise the class of profinite groups with precisely
countably infinitely many closed subgroups. We first mention three further facts
used in the proof. Firstly, the fact that if $G$ is a profinite group and $H$
is a closed subgroup of $G$ then either $H\os G$ or
$|G:H|\geqslant 2^{\aleph_0}$. 

The second fact is about properties of the $p$-adic integers,
$\mathbb{Z}_p$, for some prime $p$. As an inverse limit  of rings,
$\mathbb{Z}_p$ is a profinite ring. Every automorphism of the additive
group $\mathbb{Z}_p$ is continuous. It turns out that
$Aut(\mathbb{Z}_p)$ is topologically isomorphic to $U(\mathbb{Z}_p)$,
the group of units of $\mathbb{Z}_p$. Further, $U(\mathbb{Z}_p) \cong
\mathbb{Z}_p \times K_p$, where $K_p = C_{p-1}$ if $p$ is odd or $K_p
=C_2$ if $p=2$. (See section~3, Chapter~II of \cite{Serre}.)

The third fact concerns a small amount of cohomology we use at the end of the
proof of our characterisation theorem. Let $G$ be a profinite group and
$N\cns G$. Suppose that $N$ has a closed complement $H$, i.e., $H\cs G$,
$G=HN$ and $H\cap N=1$ and so $G=N\rtimes H$. By a \df{derivation} from
$H$ to $N$ we mean a continuous map $d\colon H\to N$ such that
$d(g_1g_2)=(d(g_1))g_2(d(g_2))$ for all $g_1,g_2\in H$. The set of derivations
from $H$ to $N$ is denoted by $\Der (H,N)$. The basic result we require is that
there is a bijection between the set of all complements by closed subgroups of
$N$ in $G$ and $\Der (H,N)$. The proof which is very similar to the abstract
case can be found in Lemma 6.2.3(a) of \cite{WilPrf}. In fact since we
actually encounter central split extensions, we are really using the first
cohomology group, $H^1(G,N)$.

\begin{thm}\label{tcount}
Let $G$ be an infinite profinite group. Then the following are equivalent.
\begin{enumerate}
\item[(i)] $G$ has precisely countably infinitely many closed subgroups.
\item[(ii)] $G$ has less than $2^{\aleph_0}$ closed subgroups.
\item[(iii)] $\zg (G)\ons G$ and there is a finite non-empty set of primes
$\{p_1,\ldots ,p_n\}$ (with $p_i\neq p_j$ for $i\neq j$) and a finite abelian
group $F$ such that $\zg (G)$ is topologically isomorphic to
$F\times\Dr_{i=1}^n\mathbb Z_{p_i}$.
\item[(iv)] $G^{\prime}$ is finite and there is a finite non-empty set of
primes $\{p_1,\ldots ,p_n\}$ (with $p_i\neq p_j$ for $i\neq j$) and a finite
abelian group $F$ such that $G/G^{\prime}$ is topologically isomorphic to
$F\times\Dr_{i=1}^n\mathbb Z_{p_i}$.
\end{enumerate}
\end{thm}

\begin{prf}
Clearly $\hbox{(i)}\Rightarrow\hbox{(ii)}$. First, we show that
$\hbox{(ii)}\Rightarrow\hbox{(iii)}$, so suppose that $G$ has less than
$2^{\aleph_0}$ closed subgroups. We shall now show that $G$ is virtually
abelian. Let $H$ be a maximal abelian subgroup of $G$ (which exists by Zorn's
lemma). Then $\overline{H}$ is also abelian, and so $H\cs G$ with
$H=\cen_G(H)$. Let 
$K=\norm_G(H)$. Then $|G:K|=|G:\norm_G(H)|=|\{H^g\mid g\in G\}|<2^{\aleph_0}$.
So by the first fact above, $K\os G$. If $H$ is finite, then
$K/H=\norm_G(H)/\cen_G(H)$ is finite, and thus $G$ is finite, a contradiction.
Hence $H$ is infinite.

Now by Lemma \ref{lcab}, $H=T\times\Dr_{i=1}^n H_i$ where $T$ is a finite
subgroup of $H$, $\{p_1,\ldots, p_n\}$ is a finite non-empty set of primes
(with $p_i\neq p_j$ for $i\neq j$) and $H_i\cong\mathbb Z_{p_i}$ for each $i$.
$T$ is a characteristic subgroup of $H$, and so $T\cns K$. Also for each $i$,
$TH_i$ is a characteristic subgroup of $H$, and so $TH_i\cns K$. Let
$C_i=\cen_K(TH_i/T)$, then $C_i\cns K$. Suppose there is an
$i$ such that $K/C_i$ is infinite. But, by our discussion above,
$K/C_i\cong\Aut (TH_i/T)\cong\Aut (\mathbb Z_{p_i})\cong\mathbb Z_{p_i}
\times K_{p_i}$, where $K_{p_i}=C_{p_i-1}$ if $p_i$ is odd or
$K_{p_i}=C_2$ if $p_i=2$. Thus there exists $L_i/C_i\cs K/C_i$ such that
$L_i/C_i\cong\mathbb Z_{p_i}$. Also $TH_i/T\cs C_i$,and
$TH_i/T\cong\mathbb Z_{p_i}$. Hence by Lemma \ref{lexsm}(iii), $L_i/T$ has at
least $2^{\aleph_0}$ closed subgroups, which is a contradiction. So for every
$i$, $C_i\ons K$. Hence $\bigcap_{i=1}^n C_i=\cen_K(H/T)\ons K$.

Clearly $\cen_K(T)\ons K$, so $L=\cen_K(H/T)\cap \cen_K(T)\ons K$. Also
$H\cs L$, so to show that $H\os G$, it suffices to show that $H\os L$. For
each $x\in L$, define $\theta_x\colon H\to T$, by $\theta_x(h)=[h,x]$. Then as
$H\cs \cen_K(T)$, $\theta_x$ is a homomorphism and is continuous 
 since $H$ is a finitely generated (and otherwise very simple) 
profinite group. Define
$\theta\colon L\to \Hom (H,T)$, by $\theta (x)=\theta_x$. Then as
$L\cs \cen_K(T)$, $\theta$ is a homomorphism. Now
$\Hom (H,T)\cong\Hom (T,T)\times\Dr_{i=1}^n\Hom (\mathbb Z_{p_i},T)$, which
is finite. So $\ker\theta\ons L$ since $L$ is strongly complete. But clearly
$\ker\theta =\cen_L(H)=H$. Hence $H\os G$.

Let $N$ be the core of $\Dr_{i=1}^n H_i$ in $G$. Then $N\ons G$ and
$N=\Dr_{i=1}^n N_i$ where $N_i\cong\mathbb Z_{p_i}$ for each $i$. To show
that $\zg (G)\ons G$ we show that $N\leqslant\zg (G)$. So suppose for a
contradiction that $N\nleqslant\zg (G)$. Then there is an $i$ such that
$G\neq\cen_G(N_i)$. Let $M=\Dr_{{j=1}\atop {j\neq i}}^n N_j$. For each $j$,
$N_j$ is characteristic in $N$ and so $N_j\cns G$. Thus $M\cns G$. Then
$G/M\neq \cen_{G/M}(N/M)$ and $N/M\cong N_i$. So we may assume that
$N_i=N\cong\mathbb Z_p$ for some prime $p$. Then $N\ons G$. Let $G$ act
on $N$ by conjugation. $N$ has a natural ring structure, and as {\sl
  per} our discussion above of $\mathbb{Z}_p$, 
$\Aut (N)\cong U(N)$, the group of units of $N$. Let
$x\in G\setminus \cen_G(N)$, and $n\in N$. There exists $u\in U(N)$ such that
$u\cdot n=n^x$ where $\cdot$ is the ring multiplication in $N$. If $n^x=n$
then $u\cdot n=n$ and $n$ is the additive identity in $N$. Thus $\cen_N(x)$ is
trivial and $N\cap\langle x\rangle =1$. Since $N\ons G$ and
$\langle x\rangle=\langle x\rangle/N\cap\langle x\rangle\cong N\langle x\rangle
/N$, $x$ is of finite order and so $\langle x\rangle\cs G$. Define
$\phi\colon N\to\{H\mid H\cs G\}$ by $\phi (g)=\langle x\rangle^g$. If
$\phi(g_1)=\phi(g_2)$ then $\langle x^{g_1}\rangle=\langle x^{g_2}\rangle$.
So there is an integer $k$ such that $x^{g_1}=(x^{g_2})^k$. Then $Nx=Nx^k$ and
as $N\cap\langle x\rangle =1$, $k=1$. Thus $g_1=g_2$. Hence $\phi$ is injective,
and as $|N|=2^{\aleph_0}$ we have a contradiction. Thus $N\leqslant\zg (G)$, and
$\zg (G)\ons G$. The structure of $\zg (G)$ is as required by Lemma \ref{lcab}.

Now for $\hbox{(iii)}\Rightarrow\hbox{(i)}$. By Lemma \ref{labvp} we have
nothing more general than the class of groups in (iii), by supposing that
$\zg (G)\ons G$, and we have an $N\os\zg (G)$ with $N=\Dr_{i=1}^n N_i$ where
$N_i\cong\mathbb Z_{p_i}$ and $\{p_1,\ldots ,p_n\}$ is a finite non-empty set
of primes with $p_i\neq p_j$ for $i\neq j$. Suppose for a contradiction that
$G$ has uncountably many closed subgroups. Let $(H_{\alpha})_{\alpha\in A}$ be
a family of closed subgroups of $G$ with $A$ uncountable, and
$H_{\alpha}\neq H_{\beta}$ for $\alpha\neq\beta$. Firstly we show that we may
assume without loss of generality that $G$ splits over $N$. Since $N\ons G$,
there are only finitely many closed subgroups of $G$ containing $N$. 
So there exists an uncountable subset $B$ of $A$ such
that $H_{\alpha}N=H_{\beta}N$ for every $\alpha$, $\beta\in B$. Let
$L=H_{\alpha}N$ for $\alpha\in B$. Also for every $\alpha\in B$,
$H_{\alpha}\cs L$, $N\os L$ and $N\os\zg (L)$. Thus without loss of generality
we may assume that $H_{\alpha}N=G$ for every $\alpha\in A$.

Now $N$ has precisely countably infinitely many closed subgroups. So
$\{H_{\alpha}\cap N\mid\alpha\in A\}$ is at most countable. Thus there exists
an uncountable subset $C$ of $A$ such that $H_{\alpha}\cap N=H_{\beta}\cap N$
for every $\alpha$, $\beta\in C$. Let $K=H_{\alpha}\cap N$ for $\alpha\in C$.
Then since $K\leqslant\zg (G)$, $K\cns G$. If $K\ons G$, then $H_{\alpha}\os G$
for every $\alpha\in C$, a contradiction as $G$ has only countably many open
subgroups. Clearly $G/K$ has uncountably many closed subgroups, $H_{\alpha}/K$
($\alpha\in C$). As $\zg (G)/K\leqslant \zg (G/K)$, $N/K\os\zg (G/K)\ons G/K$.
If $K\os N$, then $K\os G$, a contradiction. So there exists $M/K\cs N/K$ such
that $M/K$ is topologically isomorphic to $\Dr_{p\in P}\mathbb Z_p$ where $P$
is a non-empty subset of $\{p_1,\ldots ,p_n\}$. Clearly for each $\alpha\in C$,
$H_{\alpha}/K\cap M/K\leqslant (H_{\alpha}\cap N)/K=K/K$. Thus we may assume
without loss of generality that $H_{\alpha}\cap N=1$ for every $\alpha\in A$.

So now $G$ splits over $N$ and every $H_{\alpha}$ ($\alpha\in A$) is a
complement to $N$ in $G$. Let $\alpha\in A$. Then
$$
\begin{array}{lrl}
\multicolumn{3}{l}{|\{H\cs G\mid H\hbox{ is a complement to $N$ in
    $G$}\}|} \\
\phantom{XXXXXXXXXXXXX} &=& |\Der(H_{\alpha},N)| \\ 
&=& |\Hom(H_{\alpha},N)| \qquad\hbox{(as $N$ is central in $G$)} \\
&=& |\Hom(H_{\alpha},\Dr_{i=1}^n N_i)| \\
&=& |\Dr_{i=1}^n \Hom(H_{\alpha},N_i)| \\
&=& 1 \\
\end{array}
$$
This is a contradiction as every $H_{\alpha}$ is a complement to $N$ in $G$.
Hence $G$ has precisely countably infinitely many closed subgroups.

Now suppose that $G$ has precisely countably infinitely many closed
subgroups. Then we know (iii) holds. In particular $\zg (G)\ons G$. So by
Proposition \ref{pzgocomf}(i), $G^{\prime}$ is finite. Then by Lemma
\ref{lcab}, $G/G^{\prime}$ is of the required form for (iv).

Finally, suppose that (iv) holds. Then, clearly $G$ is finitely generated
and so by Proposition \ref{pzgocomf}(ii), $\zg (G)\ons G$. Now
$\zg (G)/(\zg (G)\cap G^{\prime})
\cong \zg (G)G^{\prime}/G^{\prime}\ons G/G^{\prime}$. So for each $i$,
$\h_{p_i}(\zg (G)/(\zg (G)\cap G^{\prime}))=1$. But $\zg (G)\cap G^{\prime}$
is finite and so $\h_{p_i}(\zg (G))=1$ for each $i$. Consequently it is clear
that $\zg (G)$ is of the required form for (iii), and so
$\hbox{(iv)}\Rightarrow\hbox{(iii)}$.
\end{prf}

\section{Topological Classification of Countable $\sub (G)$}

Applying the classification of countable profinite spaces of
Section~2.3, we are now able to describe $\sub (G)$ for $G$ profinite
and $\sub (G)$ 
countable. Naturally we use our description (Theorem \ref{tcount}) of such
groups. First we require a preliminary lemma whose proof is similar to
part of the proof of that theorem.

\begin{lemma}\label{lctcenfin}
Let $G$ be a profinite group with a central open subgroup $N$ such that
$N=\Dr_{i=1}^kN_i$ where $\{p_1,\ldots ,p_k\}$ is a finite non-empty set of
primes (with $p_i\neq p_j$ for $i\neq j$) and for each $i$, $N_i$ is a
procyclic pro-$p_i$ group. Let $H\cs G$. Then $\{K\cs G\mid K\cap N=H\cap N\}$
is finite.
\end{lemma}

\begin{prf}
Suppose for a contradiction that there is a profinite group $G$ satisfying
the above conditions with a closed subgroup $H$ of $G$ and an infinite family
$(K_{\lambda})_{\lambda\in\Lambda}$ of closed subgroups of $G$ such that
$K_{\lambda}\neq K_{\mu}$ for $\lambda\neq\mu$ and $K_{\lambda}\cap N=H\cap N$
for every $\lambda\in\Lambda$. We show that without loss of generality
$G$ splits over $N$. Since $N\ons G$, there are only finitely many
closed subgroups of $G$ containing $N$. So there exists an infinite subset
$\Lambda^{\prime}$ of $\Lambda$ such that $K_{\lambda}N=K_{\mu}N$ for every
$\lambda,\mu\in\Lambda^{\prime}$. Now for every $\lambda\in\Lambda^{\prime}$,
$N\leqslant\zg (K_{\lambda}N)$ and $K_{\lambda}\cs K_{\lambda}N$. So without
loss of generality we may assume that $K_{\lambda}N=G$ for every
$\lambda\in\Lambda$. Since $N$ is central in $G$, $H\cap N\cns G$. Clearly
for every $\lambda\in\Lambda$,
$(K_{\lambda}/(H\cap N))(N/(H\cap N))=G/(H\cap N)$ and
$(K_{\lambda}/(H\cap N))\cap (N/(H\cap N))$ is trivial. Also
$(N/(H\cap N))\leqslant (\zg (G)/(H\cap N))\leqslant\zg (G/(H\cap N))$.
Moreover $(N/(H\cap N))=\Dr_{i=1}^k(N_i/(H\cap N_i))$ and each
$N_i/(H\cap N_i)$ is a procyclic pro-$p_i$ group. Hence without loss of
generality we may assume that $G$ splits over $N$ and every $K_{\lambda}$
for $\lambda\in\Lambda$ is a complement to $N$ in $G$, and thus finite.

Then if $\lambda\in\Lambda$,
$$
\begin{array}{lrl}
\multicolumn{3}{l}{|\{K\cs G\mid K\hbox{ is a complement to $N$ in
    $G$}\}|} \\
\phantom{XXXXXXXXXXXXX} &=& |\Der(K_{\lambda},N)| \\ 
&=& |\Hom(K_{\lambda},N)| \qquad\hbox{(as $N$ is central in $G$)} \\
&=& |\Hom(K_{\lambda},\Dr_{i=1}^k N_i)| \\
&=& |\Dr_{i=1}^k \Hom(K_{\lambda},N_i)| \\
\end{array}
$$
But for each $i$, $\Hom(K_{\lambda},N_i)$ is finite. Thus there are only
finite many complements to $N$ in $G$. But this is a contradiction since
every $K_{\lambda}$ for $\lambda\in\Lambda$ is a complement to $N$ in $G$.
\end{prf}

\begin{thm}\label{pcountsub}
Let $G$ be a profinite group with $\sub (G)$ countably infinite. Let $k$ equal
the number of primes $p$ such that $p^{\infty}\mid o(G)$. Then there exists a
positive integer $n$ such that $\sub (G)$ is homeomorphic to
$\omega^kn +1$.
\end{thm}

\begin{prf}
Firstly note that by Theorem \ref{tcount} there are only finitely many primes
$p$ such that $p^{\infty}\mid o(G)$. By Proposition \ref{pctprs} it suffices to
show that $\hht (\sub (G))=k+1$. Now by Theorem \ref{tcount} $G$ has a central
open subgroup $N$ such that $N$ is topologically isomorphic to
$\Dr_{i=1}^k\mathbb Z_{p_i}$ for some finite non-empty set of primes
$\{p_1,\ldots, p_k\}$ (with $p_i\neq p_j$ for $i\neq j$). Now clearly by
Example \ref{esubzp} and Proposition \ref{pctprfct}, $\sub (N)$ is homeomorphic
to $(\omega +1)^k$. But $\hht ((\omega +1)^k)=k+1$ and so by Lemma
\ref{lhtopen}, $\hht (\sub (G))\geqslant k+1$. Thus it suffices to show that
for every $H\cs G$, $\hht (H,\sub (G))\leqslant k+1$.

We know that $\hht (H\cap N,\sub (N))\leqslant k+1$ and by Lemma \ref{lhtopen},
$\hht (H\cap N,\sub (G))\leqslant k+1$. Hence it suffices to show that
$\hht (H,\sub (G))\leqslant\hht (H\cap N,\sub (G))$. Now by Lemma \ref{lclbase}
there exists $U\ons G$ with $U\leqslant N$ such that if $L\in B(H\cap N,U)$
then either $L=H\cap N$ or $\hht (L,\sub (G))<ht(H\cap N,\sub (G))$. By Lemma
\ref{lctcenfin} $\{K\cs G\mid K\cap N=H\cap N\}$ is finite. Hence since
$\sub (G)$ is Hausdorff there exists $V\ons G$ with $V\leqslant U$ such that
$B(H,V)\cap\{K\cs G\mid K\cap N=H\cap N\}=\{H\}$. Now let $K\in B(H,V)$. Then
$K\in B(H,U)$ and since $U\leqslant N$, $K\cap N\in B(H\cap N,U)$. Thus
by the above choice of $U$, either $K=H\cap N$ or
$\hht (K\cap N,\sub (G))<\hht (H\cap N,\sub (G))$. But if $K=H\cap N$ then by
the above $K=H$. Suppose that $\hht (K\cap N,\sub (G))<\hht (H\cap N,\sub (G))$.
Then, by induction on $\hht (H\cap N,\sub (G))$,
$\hht (K,\sub (G))\leqslant\hht (K\cap N,\sub (G))$. So then
$\hht (K,\sub (G))<\hht (H\cap N,\sub (G))$. Hence
$\hht (H,\sub (G))\leqslant\hht (H\cap N,\sub (G))$ as required.
\end{prf}

\begin{cor}\label{cctsubprop}
Let $G$ be a pro-$p$ group with $\sub (G)$ countably infinite. Let $n$
equal the number of closed non-open subgroups of $G$. Then $\sub (G)$
is homeomorphic to $\omega n+1$.
\end{cor}

\section{Isolated Subgroups}\label{sisolsub}

In this section we characterise the profinite groups $G$ for which $\sub (G)$
is perfect (i.e. does not contain an isolated point). 
In particular we determine the profinite groups $G$ such that
$\sub (G)$ is homeomorphic to the Cantor set. We also make some
remarks on the situation for $\nsub (G)$.

\paragraph{When Is $\sub (G)$ Perfect?}

\begin{defn}
If $G$ is a profinite group and $H\cs G$ then $H$ is said to be an
\df{isolated subgroup} of $G$ if $H$ is an isolated point of $\sub (G)$.
\end{defn}

So, for example (cf. Example \ref{esubzp}), for each non-negative integer $n$,
$p^n\mathbb Z_p$ is an isolated subgroup of $\mathbb Z_p$ whereas $0$ is
not an isolated subgroup of $\mathbb Z_p$.

\begin{defn}
For $G$ a profinite group, let $\Psi(G)$ be the intersection of the maximal
(proper) open normal subgroups of $G$.

Recall that the \df{Frattini subgroup}, $\Phi(G)$, of $G$ is the
intersection of the maximal (proper) open subgroups. 
\end{defn}

We require some well-known facts about the Frattini subgroup of a profinite
group $G$ and about $\Psi(G)$. For example Lemma \ref{lfratpsi}(i) is a
profinite version of a well-known result about finite groups, cf. Lemma
11.4 of \cite{Rose}.

\begin{lemma}\label{lfratpsi}
Let $G$ be a profinite group.
\begin{enumerate}
\item[(i)] If $K\cns G$, then $K\leqslant\Phi(G)$ if and only if $H\cs G$
and $HK=G$ implies that $H=G$.
\item[(ii)] If $K\cns G$, then $K\leqslant\Psi(G)$ if and only if $H\cns G$
and $HK=G$ implies that $H=G$.
\item[(iii)] $\Phi(G)\leqslant\Psi(G)$.
\item[(iv)] If $G$ is pronilpotent then $\Phi(G)=\Psi(G)$.
\end{enumerate}
\end{lemma}

\begin{prf}
For (i), first suppose that $K\leqslant\Phi(G)$, and let $H$ be a proper
closed subgroup of $G$. Then  $H$ is contained in
a maximal open subgroup $M$ of $G$. Since $K\leqslant\Phi(G)$, $K\leqslant M$.
But now $HK\leqslant M\neq G$. Now suppose that $K\nleqslant\Phi(G)$. Now
$G$ must be non-trivial and there exists a maximal open subgroup $M$ of $G$
such that $K\nleqslant M$. But then $M$ is a proper subgroup of $MK$.
Thus $MK=G$. The argument for (ii) is basically the same.

For (iii), let $N$ be a maximal open normal subgroup of $G$. Then $N$ is
contained in a maximal open subgroup $M$ of $G$. Clearly $N=M_G$. But each
conjugate of $M$ in $G$ is a maximal open subgroup of $G$. So $N$ is an
intersection of maximal open subgroups of $G$. Thus $\Phi(G)\leqslant N$, and
so $\Phi(G)\leqslant\Psi(G)$.

For (iv), let $M$ be a maximal open subgroup of $G$. Let $N=M_G$. Then
$M/N$ is a maximal subgroup of $G/N$. But $G/N$ is a finite nilpotent
group, and so $M/N\triangleleft G/N$. Thus $M\cns G$, and so clearly
$\Phi(G)=\Psi(G)$.
\end{prf}

Note that it follows immediately from Lemma \ref{lfratpsi}(i) that if $G$
is a profinite group and $H\cs G$ with $H\Phi(G)=G$ then $H=G$ (cf.
Proposition 2.5.1(a) of \cite{WilPrf}). Of course the corresponding fact holds
for $\Psi(G)$.

We now prove some elementary facts about isolated subgroups.

\begin{lemma}\label{lifrat}
Let $G$ be a profinite group and $H\cs G$.
\mbox{\ }
\begin{enumerate}
\item[(i)] $H$ is an isolated subgroup of $G$ if and only if there exists an
open normal subgroup $N$ of $G$ such that if $K\cs G$ and $KN=HN$ then $K=H$.
\item[(ii)] If $H$ is an isolated subgroup of $G$ then $H\os G$.
\item[(iii)] If $H\os G$ then $H$ is an isolated subgroup of $G$ if and only
if there exists an open normal subgroup $N$ of $G$ with $N\leqslant H$ such
that if $K\cs G$ and $KN=H$ then $K=H$.
\item[(iv)] If $H\os G$ then $H$ is an isolated subgroup of $G$ if and only
if $\Phi(H)\ons H$.
\end{enumerate}
\end{lemma}

\begin{prf}
By Lemma \ref{lclbase}, $H$ is an isolated subgroup of $G$ if and only if
$B(H,N)=\{H\}$ for some open normal subgroup $N$ of $G$. (i) now follows
immediately.

For (ii), by (i) there is an open normal subgroup $N$ of $G$ such that if
$K\cs G$ and $KN=HN$ then $K=H$. Let $K=HN$. Then $KN=HN$, and so $HN=H$,
that is $N\leqslant H$. Thus $H\os G$.

For (iii), let $H$ be an isolated open subgroup of $G$. Then by (i) there
exists an open normal subgroup $M$ of $G$ such if $K$ is a closed subgroup
of $G$ and $KM=HM$ then $K=H$. Let $N=(H\cap M)_G$. Then $N\ons G$ and
$N\leqslant H$. Suppose that $K$ is a closed subgroup of $G$ and $KN=H$.
Then $H\leqslant KM$. Thus $HM\leqslant KM$. But $K\leqslant H$, so
$KM=HM$. Hence $K=H$. The other implication is immediate from (i).

For (iv), first suppose that $H\os G$ and $\Phi(H)\nons H$. Let $N\ons G$
with $N\leqslant H$. Then $N\nleqslant\Phi(H)$. So there exists a maximal
open subgroup $M$ of $H$ such that $N\nleqslant M$. Then $MN=H$ but
$M\neq H$. Hence $H$ is not an isolated subgroup of $G$ by (iii).

Now suppose that $H$ is an open subgroup of $G$ and $\Phi(H)\ons H$. Let
$N=\Phi(H)_G$. Then $N\ons G$ and $N\leqslant\Phi(H)$. Suppose for $K\cs G$
that $KN=HN$. Then $KN=H$, so $H\leqslant K\Phi(H)$. But also
$\Phi(H)\leqslant H$ and $K\leqslant H$, so $H=K\Phi(H)$. Thus $H=K$ by Lemma
\ref{lfratpsi}(i). Hence $H$ is an isolated subgroup of $G$ by (i).
\end{prf}

The following lemma is well known.

\begin{lemma}\label{lpfrat}
\mbox{\ }
\begin{enumerate}
\item[(i)] Let $G$ be a profinite group and $p$ be a prime. Then $p\mid o(G)$
if and only if $p\mid |G:\Phi(G)|$.
\item[(ii)] Let $G$ be a profinite group. If $G$ is finitely generated,
pronilpotent and only finitely many primes divide the order of $G$ then
$\Phi(G)\ons G$.
\end{enumerate}
\end{lemma}

\begin{prf}
For (i),  if $p\mid |G:\Phi (G)|$ then $p\mid o(G)$. So
now suppose $p\mid o(G)$ but $p\nmid |G:\Phi (G)|$. Let $P$ be a Sylow-$p$
subgroup of $\Phi (G)$. Then by a profinite version of the Frattini argument
(see Proposition 2.2.3(c) of \cite{WilPrf}), $P\cns G$. Clearly $p\nmid |G:P|$,
and so $P$ is a normal Hall subgroup of $G$. Then by the profinite version of
the Schur--Zassenhaus theorem (see Proposition 2.3.3 of \cite{WilPrf}), there
exists a closed subgroup $H$ of $G$ such that $G=P\rtimes H$. But then by Lemma
\ref{lfratpsi}(i), $H=G$, a contradiction.

For (ii), suppose $G$ is finitely generated, pronilpotent and only finitely
many primes divide the order of $G$. Then since $G$ is pronilpotent it
is topologically isomorphic to the product of its Sylow subgroups (see
Proposition~2.4.3 of \cite{WilPrf}). So there
exist primes $p_1,\ldots ,p_n$ such that $G$ is topologically isomorphic to
$\Dr_{i=1}^n G_{p_i}$ where $G_{p_i}$ is the unique Sylow-$p_i$ subgroup of
$G$. Then $G/\Phi (G)$ is topologically isomorphic to
$\Dr_{i=1}^n G_{p_i}/\Phi (G_{p_i})$ (see Lemma 20.4 of \cite{FJ} for example).
But for each $i$, $G_{p_i}$ is finitely generated and so  
$\Phi (G_{p_i})\ons G_{p_i}$ (see Proposition~1.14 \cite{FJ}). Hence
$\Phi (G)\ons G$. 
\end{prf}

\begin{thm}\label{pchperf}
Let $G$ be a profinite group. Then the following are equivalent.
\begin{enumerate}
\item[(i)] $S(G)$ is not perfect.
\item[(ii)] $G$ has an isolated open subgroup.
\item[(iii)] Every open subgroup of $G$ is isolated.
\item[(iv)] $\Phi (G)\ons G$.
\item[(v)] $G$ is finitely generated, virtually pronilpotent and only
finitely many primes divide the order of $G$.
\end{enumerate}
\end{thm}

\begin{prf}
$\hbox{(i)}\Rightarrow\hbox{(ii)}$ follows immediately from Lemma
\ref{lifrat}(ii). Clearly $\hbox{(ii)}\Rightarrow\hbox{(i)}$ and
$\hbox{(iii)}\Rightarrow\hbox{(ii)}$. It is also clear by Lemma \ref{lifrat}
that $\hbox{(iii)}\Rightarrow\hbox{(iv)}$ and that
$\hbox{(iv)}\Rightarrow\hbox{(ii)}$.

For $\hbox{(ii)}\Rightarrow\hbox{(v)}$ let $H$ be an isolated open subgroup
of $G$. Then by Lemma \ref{lifrat}(iv), $\Phi (H)\ons H$. So 
 $H$ is finitely generated and hence $G$ is finitely
generated. Also  $\Phi (H)$ is pronilpotent and
so $G$ is virtually pronilpotent. Suppose infinitely many primes divide the
order of $G$. Then infinitely many primes divide the order of $H$, and so by
Lemma \ref{lpfrat}(i), $\Phi (H)\nons H$, a contradiction.

Finally for $\hbox{(v)}\Rightarrow\hbox{(iii)}$. Suppose that $G$ is finitely
generated, virtually pronilpotent and only finitely many primes divide the
order of $G$. Let $H\os G$. Then there exists an open normal subgroup $N$
of $H$ such that $N$ is pronilpotent. Now $N$ is finitely generated and only
finitely many primes divide the order of $N$, so by Lemma \ref{lpfrat}(ii),
$\Phi (N)\ons N$. Thus $\Phi (N)\os H$. But $\Phi (N)\leqslant \Phi (H)$. So
$\Phi (H)\ons H$. Hence $H$ is an isolated subgroup of $G$ by Lemma
\ref{lifrat}(iv).
\end{prf}

\begin{cor}\label{csubcant}
Let $G$ be a profinite group. Then $\sub (G)$ is homeomorphic to the Cantor
set if and only if $G$ is countably based and not (finitely generated,
virtually pronilpotent and only finitely many primes divide the order of $G$).
\end{cor}

\begin{prf}
This follows immediately from Proposition \ref{pchperf}, Corollary
\ref{csubwt} and Proposition \ref{pchcant}.
\end{prf}

\paragraph{\sl Remark} There is another way of associating a profinite
space to a profinite group $G$. 
This is in terms of the Burnside algebra of $G$. Pierce has shown (see Theorem
6.6 of \cite{Pierba}) that the space obtained, $\Sigma (G)$ is not perfect
if and only if $G$ is finitely generated, virtually pronilpotent and only
finitely many primes divide the order of $G$. So by Proposition \ref{pchperf},
this is precisely when $\sub (G)$ is not perfect. We do not know how, in
general $\Sigma (G)$ is related to $\sub (G)$; though $\Sigma (G)$ can be
constructed from the lattice of open subgroups of $G$.

\paragraph{When Is $\nsub(G)$ Perfect?}
We now consider what can be said about when $\nsub (G)$ is perfect for
$G$ a profinite group. It is clear that by the same arguments we have a
direct analogue of Lemma \ref{lifrat} using $\Psi (G)$. In particular we have
the following.

\begin{lemma}\label{lipsi}
Let $G$ be a profinite group and $H\cns G$.
\begin{enumerate}
\item[(i)] $H$ is isolated in $\nsub (G)$ if and only if there exists an
open normal subgroup $N$ of $G$ such that if $K\cns G$ and $KN=HN$ then
$K=H$.
\item[(ii)] If $H$ is isolated in $\nsub (G)$ then $H\ons G$.
\item[(iii)] If $H\ons G$ then $H$ is isolated in $\nsub (G)$ if and only if
$\Psi(H)\ons H$.
\end{enumerate}
\end{lemma}

\begin{example}\label{big_is_n}
 If $\kappa$ is any infinite cardinal then a profinite group  with a
 unique maximal open normal subgroup  of weight
$\kappa$ can be constructed as a split extension of the Cartesian product
of $\kappa$ copies of a finite simple group, by $A_5$. 

For such a group $G$,
clearly $\Psi (G)\ons G$, and thus by Lemma \ref{lipsi}, $G$ is isolated in
$\nsub (G)$. 
\end{example}
Hence there exist infinite profinite groups of arbitrary weight
with $\nsub (G)$ not perfect. 
So we have no simple analogue of Theorem~\ref{pchperf} for the space
of normal subgroups. We can say the following. 

\begin{propn}\label{psubnperf}
Let $G$ be a profinite group. If $\nsub (G)$ is perfect then $\sub (G)$
is perfect. If $G$ is pronilpotent and $\sub (G)$ is perfect then $\nsub (G)$
is perfect.
\end{propn}

\begin{prf}
Suppose that $\sub (G)$ is not perfect. Then by Proposition \ref{pchperf},
$\Phi(G)\ons G$. Thus by Lemma \ref{lfratpsi}(iii), $\Psi(G)\ons G$. So,
by Lemma \ref{lipsi}(iii), $G$ is isolated in $\nsub (G)$ and $\nsub (G)$ is
not perfect.

Now suppose that $G$ is pronilpotent. If $H\cs G$ then by Lemma
\ref{lfratpsi}(iv), $\Phi(H)=\Psi(H)$. It is now clear from the proof of
Proposition \ref{pchperf} that $\sub (G)$ is perfect if and only if
$\nsub (G)$ is perfect.
\end{prf}

\section{Profinite Groups of Large Weight}\label{sspcancub}

In this section we examine the situation when our profinite group,
$G$, is not countably based, $w(G) > \aleph_0$. 
By Theorem~\ref{pchperf} we know $\sub (G)$ is perfect. Since for
countably based $G$ 
with perfect $\sub (G)$,
the space of subgroups is homeomorphic to $\{0,1\}^{\aleph_0}$, the
natural conjecture -- solving simultaneously the problem of counting
and topologically classifying $\sub (G)$ -- is that $\sub(G)$ is
homeomorphic to $\{0,1\}^{w(G)}$. So we start this section by
introducing some topological properties (Dugundji compactness and
character homogeneity) which  characterise 
the Cantor cube $\{0,1\}^{\kappa}$.
Then we note that $\sub (G)$ always has a property
($\kappa$-metrizability) just a little
weaker than Dugundji compactness, and is always character homogeneous.
We deduce that the conjecture is correct when $w(G)=\aleph_1$, but is
false in general. Nevertheless from our topological considerations we
have enough to deduce that $\sub(G)$ and $\nsub(G)$ both have
cardinality equal to $2^{w(G)}$. We conclude with a sketch of a more
algebraic proof of the latter fact for $\sub (G)$. However we know of
no algebraic proof that $|\nsub (G)|=2^{w(G)}$.

\paragraph{Topological Properties of Cantor Cubes} 
\begin{defns}
Let $X$ be a topological space.

$X$ is said to be \df{Dugundji compact} if $X$ is compact and Hausdorff
and whenever $Y$ is a profinite space, $Z$ is a closed subspace of $Y$,
and $f\colon Z\to X$ is a continuous map, then $f$ can be extended to
a continuous map $f^{\prime}\colon Y\to X$.

Now let $\mathcal{RC}(X)$ denote the family of all regularly closed
subsets of $X$; that is
$\mathcal{RC}(X)=\{F\subseteq X\mid F=\overline{\Int F}\}$.

A $\kappa$-metric on $X$ is a non-negative function
$\rho\colon X\times\mathcal{RC}(X)\to\mathbb R$ satisfying the following
four conditions:
\begin{enumerate}
\item[(i)] $\rho (x,F)=0$ if and only if $x\in F$;
\item[(ii)] if $x\in X,F,F^{\prime}\in\mathcal{RC}(X)$ and
$F\subseteq F^{\prime}$, then $\rho (x,F^{\prime})\leqslant\rho (x,F)$;
\item[(iii)] given $F\in\mathcal{RC}(X)$, the function
$\rho(\cdot,F)\colon X\to\mathbb R$ is continuous with respect to the
first argument;
\item[(iv)] if $(F_{\alpha})_{\alpha<\tau}$ is an increasing transfinite
sequence of regularly closed sets in $X$, and $x\in X$, then
$\rho (x,\overline{\bigcup_{\alpha<\tau}F_{\alpha}})=
\inf_{\alpha<\tau}\rho (x,F_{\alpha})$.
\end{enumerate}

$X$ is said to be \df{$\kappa$-metrizable} if $X$ admits a $\kappa$-metric.

$X$ is said to be \df{$\kappa$-adic} if there is a $\kappa$-metrizable
space $Y$ and a continuous surjection $f\colon Y\to X$.
\end{defns}

The definition of Dugundji compact spaces we have given is not the original
definition due to Pe\l czy\'nski. Our definition asserts that a Dugundji
compact space is a compact Hausdorff space with the property of being an
absolute extensor in dimension $0$ (AE($0$)). The equivalence of these
definitions is due to Haydon, \cite{Hayd}. Note that a profinite space is
Dugundji compact if and only if it is an injective object in the category
of profinite spaces. The class of $\kappa$-metrizable spaces and the class
of $\kappa$-adic spaces are due to Shchepin, \cite{Shch1}, \cite{Shch2}.
\S 7 of Shakhmatov, \cite{Shak} is a survey of these three class of spaces.

We also require the following simple terms:
\begin{defns}
Let $X$ be a topological space.

A \df{$\pi$-basis for the open neighbourhoods} of $x\in X$ is a collection
$\mathcal B (x)$ of open sets in $X$ such that every open neighbourhood of
$x$ contains an element of $\mathcal B (x)$.

The \df{$\pi$-character} of $X$ at $x\in X$ is defined to be
$\pi\chi(x,X)=\min\{|\mathcal B(x)|\mid\mathcal B(x)$ is a $\pi$-basis for
the open neighbourhoods of $x\}$.

$X$ is said to be \df{character homogeneous} if $\chi (x,X)=\chi (X)$ for
every $x\in X$.
\end{defns}
In the next theorem we collect together all the topological results we
require for this section.
\begin{thm}\label{ttfshchshap}
\mbox{\ }
\begin{enumerate}
\item[] (Shchepin)
\item[(i)] A compact Hausdorff space $X$ if $\kappa$-metrizable if and only
if $X$ is homeomorphic to the inverse limit of an inverse system
$(X_i,p_{ij})$ of compact metrizable spaces indexed by a directed set $I$
such that every $p_{ij}$ is an open continuous surjective map and every
countable chain of elements of $I$ has its least upper bound in $I$.
\item[(ii)] Let $X$ be a compact Hausdorff space. If $X$ is Dugundji
compact then $X$ is $\kappa$-metrizable. If $X$ is $\kappa$-metrizable and
$w(X)=\aleph_1$, then $X$ is Dugundji compact.
\item[(iii)] Let $\kappa$ be an infinite cardinal and $X$ be a profinite
space with $w(X)=\kappa$. Then $X$ is homeomorphic to $\{0,1\}^{\kappa}$ if
and only if $X$ is Dugundji compact and character homogeneous.
\item[(iv)] If $X$ is a compact Hausdorff $\kappa$-adic space then
$\chi (X)=w(X)$.
\item[(v)] If $X$ is a compact Hausdorff $\kappa$-metrizable space and
$x\in X$ then $\pi\chi (x,X)=\chi (x,X)$.

(Shapirovski\u\i)
\item[(vi)] Let $X$ be a compact Hausdorff $\kappa$-metrizable space.
Then there exists a continuous surjective map $f\colon X\to [0,1]^{w(X)}$
(where $[0,1]^{w(X)}$ is the Hilbert cube of weight $w(X)$) if and only if
$\pi\chi (x,X)=w(X)$ for some $x\in X$.
\end{enumerate}
\end{thm}

\begin{prf}
For (i), see \cite{Shch1} and Theorem 21 of \cite{Shch3}. For (ii), see
Corollary 1 of Theorem 4 and Theorem 5 of \cite{Shch2}. For (iii), see
Theorem 9 of \cite{Shch1}. For (iv), see the Corollary to Theorem 11 of
\cite{Shch2}. For (v), see Theorem 8 of \cite{Shch2}. For (vi), see
\cite{Shap}, and also Theorem 7.21(vi) of \cite{Shak}.
\end{prf}

By (v), it might seem unnecessary to consider $\pi$-characters in (vi).
Shapirovski\u{\i} proves (vi) for a wider class of spaces, where characters
and $\pi$-characters may not coincide. We have introduced $\pi$-characters
since in our application of (vi), we can avoid the use of (v); see the
remark before Lemma \ref{lngpichwt}.

\paragraph{Topological Properties of $\sub(G)$ and $\nsub(G)$}

\begin{thm}[Fisher and Gartside, \cite{FG2}]\label{psubkm}
Let $G$ be a profinite group. Then $\sub (G)$ and $\nsub (G)$ are
$\kappa$-metrizable.
\end{thm}

\begin{lemma}\label{lchargwtg}
Let $G$ be a profinite group with $\sub (G)$ perfect. Then
$\chi (G,\sub (G))=w(G)$. If $G$ is also pronilpotent then
$\chi (G,\nsub (G))=w(G)$.
\end{lemma}

\begin{prf}
$\chi (G,\sub (G))\leqslant w(\sub (G))=w(G)$, Corollary \ref{csubwt}. By
Lemma \ref{lclbase} there exists a family
$(N_{\lambda})_{\lambda\in\Lambda}$ of open normal subgroups of $G$ with
$|\Lambda|=\chi (G,\sub (G))$ and such that
$(B(G,N_{\lambda}))_{\lambda\in\Lambda}$ is a base for the open neighbourhoods
of $G$ in $\sub (G)$. As $\sub(G)$ is Hausdorff,
$\bigcap_{\lambda\in\Lambda}B(G,N_{\lambda})=\{G\}$. Hence
if $H\cs G$ and $HN_{\lambda}=G$ for every $\lambda\in\Lambda$ then $H=G$.
Let $K=\bigcap_{\lambda\in\Lambda}N_{\lambda}$. Then if $H\cs G$ with
$HK=G$ then $HN_{\lambda}=G$ for every $\lambda\in\lambda$ and so $H=G$. So
by Lemma \ref{lfratpsi}(i) $K\leqslant\Phi (G)$. As $\sub (G)$ is perfect,
by Proposition \ref{pchperf}, $\Phi (G)\nons G$. So  $w(G/K)\geqslant
w(G/\Phi (G))=w(G)$ (the last equality follows from the fact that
density and weight coincide in a profinite group, and properties of the
Frattini subgroup, see Proposition~5.2.3(a) \cite{WilPrf}). But $G/K$ embeds in
$\Cr_{\lambda\in\Lambda}G/N_{\lambda}$ and,
$w(\Cr_{\lambda\in\Lambda}G/N_{\lambda})=|\Lambda|$. Hence
$w(G)=|\Lambda|=\chi (G,\sub (G))$.

Now suppose that $G$ is also pronilpotent. Then by Lemma \ref{lfratpsi}(iii)
$\Psi (G)=\Phi (G)$. So we may use essentially the same argument above but
using (ii) of Lemma \ref{lfratpsi} and choosing
$(N_{\lambda})_{\lambda\in\Lambda}$ such that
$(B(H,N_{\lambda})\cap\nsub (G))_{\lambda\in\Lambda}$ is a base for the
open neighbourhoods of $G$ in $\nsub (G)$.
\end{prf}

The next lemma follows immediately from Lemma \ref{lsgchwt}, Proposition
\ref{psubkm} and Theorem \ref{ttfshchshap}(v). We give a direct proof, since
it is straightforward and we can then avoid using Theorem \ref{ttfshchshap}(v).
The proof is very similar to the proof of Lemma \ref{lsgchwt}.

\begin{lemma}\label{lngpichwt}
Let $G$ be an infinite profinite group. Then $\pi\chi (1,\nsub (G))=w(G)$.
\end{lemma}

\begin{prf}
$\pi\chi (1,\nsub (G))\leqslant\chi (1,\nsub (G))=w(G)$, by Lemma
\ref{lsgchwt}. Now there exists a family $(N_{\lambda})_{\lambda\in\Lambda}$
of open normal subgroups of $G$ and a family
$(H_{\lambda})_{\lambda\in\Lambda}$ of closed normal subgroups
of $G$ such that
$(B(H_{\lambda},N_{\lambda})\cap\nsub (G))_{\lambda\in\Lambda}$ is a
$\pi$-base for the open neighbourhoods of $1$ in $\nsub (G)$ and
$|\Lambda|=\pi\chi (1,\nsub (G))$.
Let $N\ons G$. Then $B(1,N)\supseteq B(H_{\lambda},N_{\lambda})$ for some
$\lambda\in\Lambda$. Clearly
$H_{\lambda}N_{\lambda}\in B(H_{\lambda},N_{\lambda})$ and so
$H_{\lambda}N_{\lambda}\in B(1,N)$; that is
$H_{\lambda}N_{\lambda}\leqslant N$. Thus $N_{\lambda}\leqslant N$ and so
$(N_{\lambda})_{\lambda\in\Lambda}$ is a base for the open neighbourhoods
of $1$ in $G$. Hence
$w(G)=\chi (1,G)\leqslant |\Lambda|=\pi\chi (1,\nsub (G))$, and so
$\pi\chi (1,\nsub (G))=w(G)$ as required.
\end{prf}

\begin{propn}\label{psgchihom}
Let $G$ be a profinite group with $\sub (G)$ perfect. Then $\sub (G)$ is
character homogeneous.
\end{propn}

\begin{prf}
For every $H\cs G$, $\chi (H,\sub (G))\leqslant w(\sub (G))=w(G)$ by Corollary
\ref{csubwt}. Suppose for a contradiction that there exists $H\cs G$ with
$\chi (H,\sub (G))<w(G)$. Then $\chi (H,H^G)<w(G)$ where $H^G$ is the space
of conjugates in $G$ of $H$. Now, clearly $H^G$ is homogeneous, and is
isomorphic, as a $G$-space, to the coset space $G/\norm_G(H)$.
Thus $\chi (G/\norm_G(H))<w(G)$. But, by Proposition \ref{psubkm}, $\sub (G)$
is $\kappa$-metrizable, and so $G/\norm_G(H)$ is $\kappa$-adic. Hence, by
Theorem \ref{ttfshchshap}(iv), $w(G/\norm_G(H)<w(G)$. Thus $w(\norm_G(H))=w(G)$.

Now $\sub (\norm_G(H)/H)$ embeds in $\sub (\norm_G(H))$, the trivial subgroup
of $\norm_G(H)/H$ being mapped to $H$. So by Lemma \ref{lsgchwt}
$w(\norm_G(H)/H))\leqslant\chi (H,\sub(\norm_G(H))\leqslant\chi (H,\sub (G))$.
Now,  $w(G)=w(\norm_G(H))=w(\norm_G(H)/H)w(H)$. If
$\sub (H)$ is not perfect, then by Proposition \ref{pchperf}, $H$ is countably
based. Hence $w(G)=w(\norm_G(H)/H)\leqslant\chi (H,\sub (G))$, a contradiction.
So $\sub (H)$ must be perfect. But, then by Lemma \ref{lchargwtg},
$w(H)=\chi (H,\sub (H))\leqslant\chi (H,\sub (G))$. Hence
$w(G)=w(\norm_G(H)/H)w(H)\leqslant\chi (H,\sub (G))$, a contradiction.
\end{prf}

\paragraph{Counting and Classifying Large Subgroup Spaces}

\begin{thm}\label{tsubpgwta1}
Let $G$ be a profinite group with $w(G)=\aleph_1$. Then $\sub (G)$ is
homeomorphic to $\{0,1\}^{\aleph_1}$.
\end{thm}

\begin{prf}
By Proposition \ref{psubkm} and Theorem \ref{ttfshchshap}(ii), $\sub (G)$ is
Dugundji compact. By Proposition \ref{pchperf}, $\sub (G)$ is perfect and so
by Proposition \ref{psgchihom}, $\sub (G)$ is character homogeneous. The result
now follows from Theorem \ref{ttfshchshap}(iii).
\end{prf}

 Haydon (private communication) 
has shown that $\sub (C_2^{\aleph_2})$ is not Dugundji
compact. In particular $\sub (C_2^{\aleph_2})$ is not homeomorphic to
$\{0,1\}^{\aleph_2}$.

\begin{thm}\label{tncbpgcns}
Let $G$ be a non-countably based profinite group. Then
$|\sub(G)|=|\nsub (G)|=2^{w(G)}$. 
\end{thm}

\begin{prf}
We know that $|\nsub (G)|\leqslant |\sub(G)| \leqslant 2^{w(G)}$. By
Proposition \ref{psubkm}, 
$\nsub (G)$ is $\kappa$-metrizable. Also $\pi\chi (1,\nsub (G))=w(G)$ by
Lemma \ref{lngpichwt} (or by Lemma \ref{lsgchwt} and Theorem
\ref{ttfshchshap}(v)). Hence, by Theorem \ref{ttfshchshap}(vi), $\nsub (G)$
maps continuously onto $[0,1]^{w(G)}$. The result then follows.
\end{prf}

\begin{cor}\label{ccount}
For any profinite group $G$:

1) $\sub (G)$ is either countable or of size $2^{w(G)}$. and

2) $\nsub (G)$ is either countable or of size $2^{w(G)}$.
\end{cor}
\begin{prf} Let $G$ be any profinite group. 
If $G$  is not 
countably based, then by the previous theorem, $|\sub(G)|=2^{w(G)} \ge
2^{\aleph_0}$, as required.  Hence if
 $G$ has $< 2^{\aleph_0}$ many closed subgroups, then it
must be countably based. In which case  $\sub (G)$ is a compact countably based
space, and  compact countably based spaces can only be countable or of
size $2^{\aleph_0}=2^{w(G)}$, again as required.

The proof for the space of closed normal subgroups is identical.
\end{prf}

\paragraph{An Algebraic Proof that $|\sub(G)| = 2^{w(G)}$}

\begin{propn}\label{pctsgncb}
Let $G$ be a non-countably based profinite group. Then $G$ has
precisely $2^{w(G)}$ closed subgroups. In fact $G$ has precisely
$2^{w(G)}$ procyclic subgroups.
\end{propn}

\begin{prf}
It suffices to show that $G$ has $2^{w(G)}$ procyclic
subgroups. Let $(X_{\lambda})_{\lambda\in\Lambda}$ be a partition of $G$ such
that for $x,y\in G$, $\overline{\langle x\rangle}=\overline{\langle y\rangle}$
if and only if $x,y\in X_{\lambda}$ for some $\lambda\in\Lambda$. As
$G$ is uncountable, $|G|=2^{w(G)}$. Thus
$2^{w(G)}=\sum_{\lambda\in\Lambda}|X_{\lambda}|=
|\Lambda|.\sup_{\lambda\in\Lambda}|X_{\lambda}|$. Each procyclic subgroup
is countably based and thus has cardinality at
most $2^{\aleph_0}$. Hence $|X_{\lambda}|\leqslant 2^{\aleph_0}$ for each
$\lambda\in\Lambda$. Thus
$\sup_{\lambda\in\Lambda}|X_{\lambda}|\leqslant 2^{\aleph_0}$. Consequently,
$|\Lambda|\leqslant 2^{w(G)}\leqslant |\Lambda|.2^{\aleph_0}$. But
$w(G)>\aleph_0$ and so $|\Lambda|=2^{w(G)}$ as required.
\end{prf}

\bibliographystyle{plain}

\bibliography{sgbib}



\end{document}